\documentclass{amsart}
\usepackage[pdftex]{color,graphicx}
\usepackage{wrapfig}
\usepackage[numbers,square]{natbib}

\def\R{\mathbb{R}}
\def\N{\mathbb{N}}

\title[Geometry of Uniformly Rotating Stars]{A Remark on the Geometry of Uniformly Rotating Stars}
\author[S. Chanillo]{Sagun Chanillo}
\address{Department of Mathematics, Hill Center 
Rutgers, The State University Of New Jersey, 
110 Frelinghuysen Rd., Piscataway, NJ 08854-8019}
\email{chanillo@math.rutgers.edu}
\author[G.S. Weiss]{Georg S. Weiss}
\address{Department of Mathematics, Heinrich Heine University, 40225 D\"usseldorf}
\email{weiss@math.uni-duesseldorf.de}
\thanks{$2000$ {\it Mathematics Subject Classification.\/} Primary
35R35, Secondary 35J60.}
\thanks{S. Chanillo was supported in part by NSF grant DMS-0855541. We wish to thank the hospitality of the Mathematisches Forschungsinstitut Oberwolfach where we carried out some of the work 
resulting in this paper.}
\date{}
\theoremstyle{plain}
\newtheorem*{thma}{Theorem A}

 \theoremstyle{definition}
 \theoremstyle{example}

\theoremstyle{definition}

\numberwithin{equation}{section}
\begin{document}
\maketitle
\begin{abstract}
In this paper we classify the free boundary associated to equilibrium configurations of compressible, self-gravitating fluid masses, rotating with constant angular velocity. The equilibrium configurations are all critical points of an associated functional and not necessarily minimizers. Our methods also apply to alternative models in the literature where the angular momentum per unit mass is prescribed. The typical physical model our results apply to is that of uniformly rotating white dwarf stars.
\end{abstract}
\section{Introduction}
In this paper we study the free boundary associated to rotating star models of white dwarf stars
with prescribed constant angular velocity. Thus we are considering figures
of equilibrium for compressible, self-gravitating fluid masses.

There has been a tremendous amount of work on incompressible, self-gravitating fluid masses
rotating with prescribed constant angular velocity since the primary investigations
by Newton. Various mathematicians like MacLaurin, Jacobi, Dirichlet, Riemann, Poincar\'e,
H. Cartan and Chandrasekhar made significant contributions to the field, studying bifurcation
sequences and analyzing the stability of various equilibrium shapes. A historical account
and details of these investigations may be found in Chandrasekhar's treatise \cite{C2}
and Tassoul's book \cite{tass}.

In the compressible case for the model with prescribed constant angular velocity $\omega>0$
(cf. \cite{yanyan}), we consider the functional
\begin{equation}\label{J}
J(\rho)=\int_{\R^3} A(\rho) \> d\xi- {1\over 2} \int_{\R^3} \omega^2 r^2 \rho  \> d\xi- 
\int_{\R^3} \rho B\rho\> d\xi,\end{equation}
where $\xi=(\xi_1,\xi_2,\xi_3), r=\sqrt{\xi_1^2+\xi_2^2}, \rho\ge 0, \rho\in L^1(\R^3)\cap L^\infty(\R^3)$ and
$$(B\rho)(\xi)= \int_{\R^3} \frac{\rho(\eta)}{|\xi-\eta|} \> d\eta.$$
Moreover
we impose the
constraint
\begin{equation}\label{constraint}
\int_{\R^3}\rho=1.
\end{equation}
As $\rho(\xi)$ represents the density of the stellar material, 
(\ref{constraint}) means that the mass of the star is prescribed.
Later in our paper we will assume in addition that the density $\rho(\xi)$
is axisymmetric, i.e. $$\rho(\xi)=\rho(r,z)\textrm{ where } z(\xi)=\xi_3.$$
This assumption means that the star is rotating about the $\xi_3$-axis.
The function $A(\rho)$ is the pressure and thus represents the equation of state
of the stellar material. We assume that $A(\rho)\in C^1([0,+\infty))$ with $A$ strictly convex
so that $A'(\rho)$ is invertible. Further conditions on $A(\rho)$ will be stipulated below.
The first term in (\ref{J}) represents then the {\em internal energy} of the star, the second
term the {\em rotational kinetic energy} and the last term the {\em gravitational potential energy}.

In \cite{yanyan} the existence of minimizers of $J$ under the constraint
(\ref{constraint}) has been obtained. \cite{libao} contains further results
for this model of prescribed angular velocity.
In \cite{CL} support estimates for critical points of (\ref{J}) under the constraint
(\ref{constraint}) have been shown.
In particular, \cite[Theorem 1]{CL} states that for $\omega\ge \omega_0>0$,
the support of $\rho$ is contained in a ball $B_\sigma(0,0,\xi_3)$
for some $\xi_3$, where $\sigma=\sigma(\omega_0)$.
It follows that $$0\le B\rho\le C \textrm{ in } \R^3.$$
Furthermore, \cite[Theorem 2]{CL} shows that the number of connected components
of the set $\{ \rho>0\}$ is finite for any critical point $\rho$.

Critical points of $J$ with the constraint (\ref{constraint}) are
according to \cite[(0.6)]{CL} characterized by the problem
\begin{equation}\label{eq}
\rho\textrm{ is continuous and nonnegative and } 
A'(\rho)-{1\over 2} \omega^2 r^2 - B\rho = \lambda(\omega)
\textrm{ in } \{ \rho>0\},
\end{equation}
where $\lambda(\omega)$ is a Lagrange multiplier arising from the
constraint (\ref{constraint}). The focus in this paper
is to study the free boundary $\partial \{ \rho>0\}$ arising from
(\ref{eq}).

There is another model of rotating stars which has been studied in the literature,
where the angular momentum per unit mass is prescribed. Existence of 
minimizers for this alternative model has been obtained in \cite{beals},
and the study of critical points has been carried out in \cite{luo}.
Caffarelli-Friedman investigated in \cite{CF} the free boundary of minimizers
for this alternative model. As Caffarelli-Friedman deal with minimizers,
they are able to apply rearrangement methods to their functional
to obtain solutions that are increasing in one direction which simplifies
the analysis as well as the result. Unfortunately this technique does not
work for critical points in either model and creates a difficulty for our
analysis. Let us remark that the proofs presented in this paper
for critical points of $J$ with the constraint (\ref{constraint}) work
equally well for the study of the free boundary of critical points
in the model in \cite{CF}.

The principal difficulty we encounter in our classification of singularities
of the free boundary is that the nonlinearity is {\em not} an increasing function
of the solution, so that various methods stemming from the well-known
obstacle problem do not apply. Neither does the monotonicity formula
derived in \cite{ACF}. Let us also mention that our problem cannot be
transformed into the type of problems studied in \cite{partial}, so we cannot
use those results either. Another difficulty is that our equation is inhomogeneous.
In particular, the leading order term on the right-hand side is not of the form
$f(u)$. This ---together with a higher order degeneracy--- distinguishes
the present problem also from the recently researched ``unstable obstacle problem'' (see \cite{MW}, \cite{AW}, \cite{ASW} and \cite{ASW2}).

Last, let us point out that ---due to the fact that the free boundary $\partial \{ \rho>0\}$
does not necessarily coincide completely with the free boundary of the PDE problem
obtained by transformation--- we obtain in our classification of singularities
several cases later called ``pseudo cases.'' We suggest that in the case of
minimizers, rearrangement techniques similar to those used in \cite{CF}
may be used to show that solutions are decreasing in a certain direction, thus
ruling out the pseudo cases.

Setting
$u = {1\over 2} \omega^2 r^2 + B\rho + \lambda(\omega)$,
we obtain in the set $\overline{\{ \rho >0\}}$
that $u=A'(\rho)=\Phi^{-1}(\rho),$ where $\Phi:[0,+\infty)\to \R$ 
is an increasing function 
satisfying according to the
asymptotics 
$$
A(\rho)= c_1 \rho^{5\over 3} +o(\rho^{5\over 3})
\textrm{ as } \rho\to 0, A(\rho)= c_2 \rho^{4\over 3}+o(\rho^{4\over 3})
\textrm{ as } \rho\to +\infty$$
(where $c_1,c_2$ are positive constants) from Chandrasekhar's book \cite[Chapter 10]{C} and \cite[(0.2)]{CL} the asymptotic relations
\begin{align}
A'(\rho)= {5\over 3} c_1 \rho^{2\over 3} + o(\rho^{2\over 3})
\textrm{ as } \rho\to 0,& A'(\rho)= {4\over 3}c_2 \rho^{1\over 3} + o(\rho^{1\over 3})
\textrm{ as } \rho\to 0,
\\
\lim_{z\to 0+} z^{-3/2} \Phi(&z) = c >0.\end{align}
It follows (cf. \cite[(3.3)]{CL}) that
$$\Delta u = 3\omega^2 -4\pi \rho = 3\omega^2 -4\pi \Phi(u)
\textrm{ in } \{ \rho >0\}$$
and 
$$\Delta u = 3\omega^2\textrm{ in } \{ \rho = 0\}.$$
Note that as $u=A'(\rho)=\Phi^{-1}(\rho)$ is only valid in the set $\overline{\{ \rho >0\}}$,
we obtain $\{ \rho>0\} \subset \{ u>0\}$ but not necessarily
the opposite inclusion. It is however true that 
$\partial\{ \rho>0\} =\partial\{ u>0\}$ and that
if $\rho(x^0)>0$ then the connected component
of $\{ \rho>0\}$ containing $x^0$ coincides with the connected component
of $\{ u>0\}$ containing $x^0$.

Normalizing the equation as well as $\Phi$ we obtain the free boundary problem
\begin{equation}\label{solution}
\Delta u = 1 -\Phi(u)\chi_\Omega \textrm{ with an open set }\Omega \textrm{ satisfying }\chi_\Omega\le \chi_{\{u>0\}};
\end{equation}   
the equation is to be understood in the sense of distributions.

\begin{thma}
For each solution $u$ of problem (\ref{solution}) the following holds:
Apart from the singular set
$S=\{ u=0\}\cap \{ \nabla u=0\}$ the level set $\{ u=0\}$ and the boundary
$\partial \Omega$ are locally $C^{2,\alpha}$-curves.
The singular set contains in each bounded subset of $\R^2$ at most finitely
many singular points $x^0$ with the following possible asymptotics (after rotation
of the coordinate system):

1. $u(x^0+tx)/t^2\to {1\over 2}x_1^2$ as $t\to 0$, and $\{ u>0\}$ and $\Omega$ satisfy
the asymptotics demonstrated in Figure \ref{figi1}.
\begin{figure}
\begin{center}
\includegraphics[width=10cm]{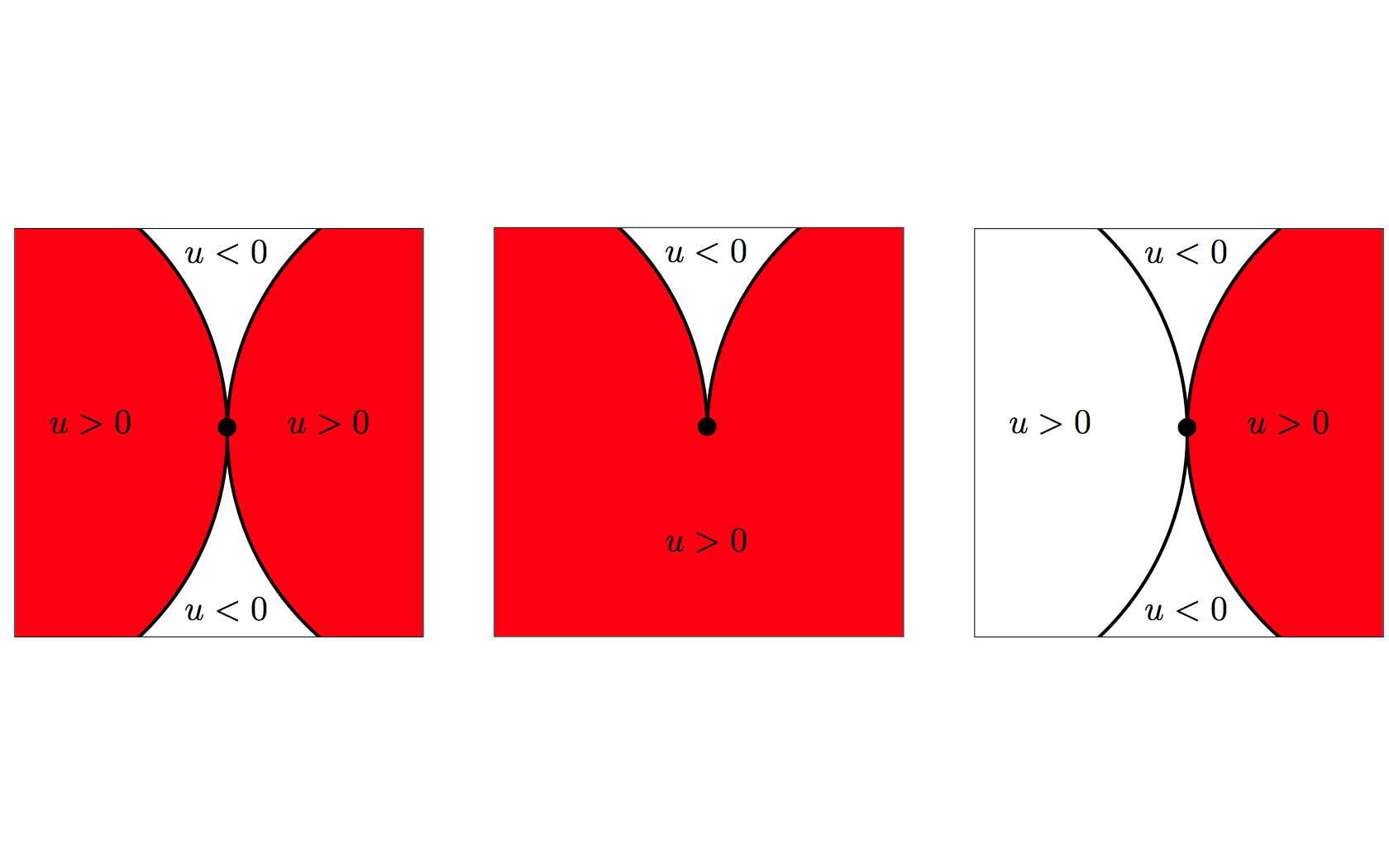}
\end{center}
\caption{Cusps (with the set $\Omega$ painted red)}\label{figi1}
\end{figure}

2. There is $\lambda\in (-\infty,0)\cup(1,+\infty)$ such that 
$u(x^0+tx)/t^2\to {1\over 2}\left(\lambda x_1^2 + (1-\lambda)x_2^2\right)$ 
as $t\to 0$. We may assume that $\lambda>1$ in which case $\{ u>0\}$ and $\Omega$ satisfy
the asymptotics demonstrated in Figure \ref{figi2}.
\begin{figure}
\begin{center}
\includegraphics[width=8cm]{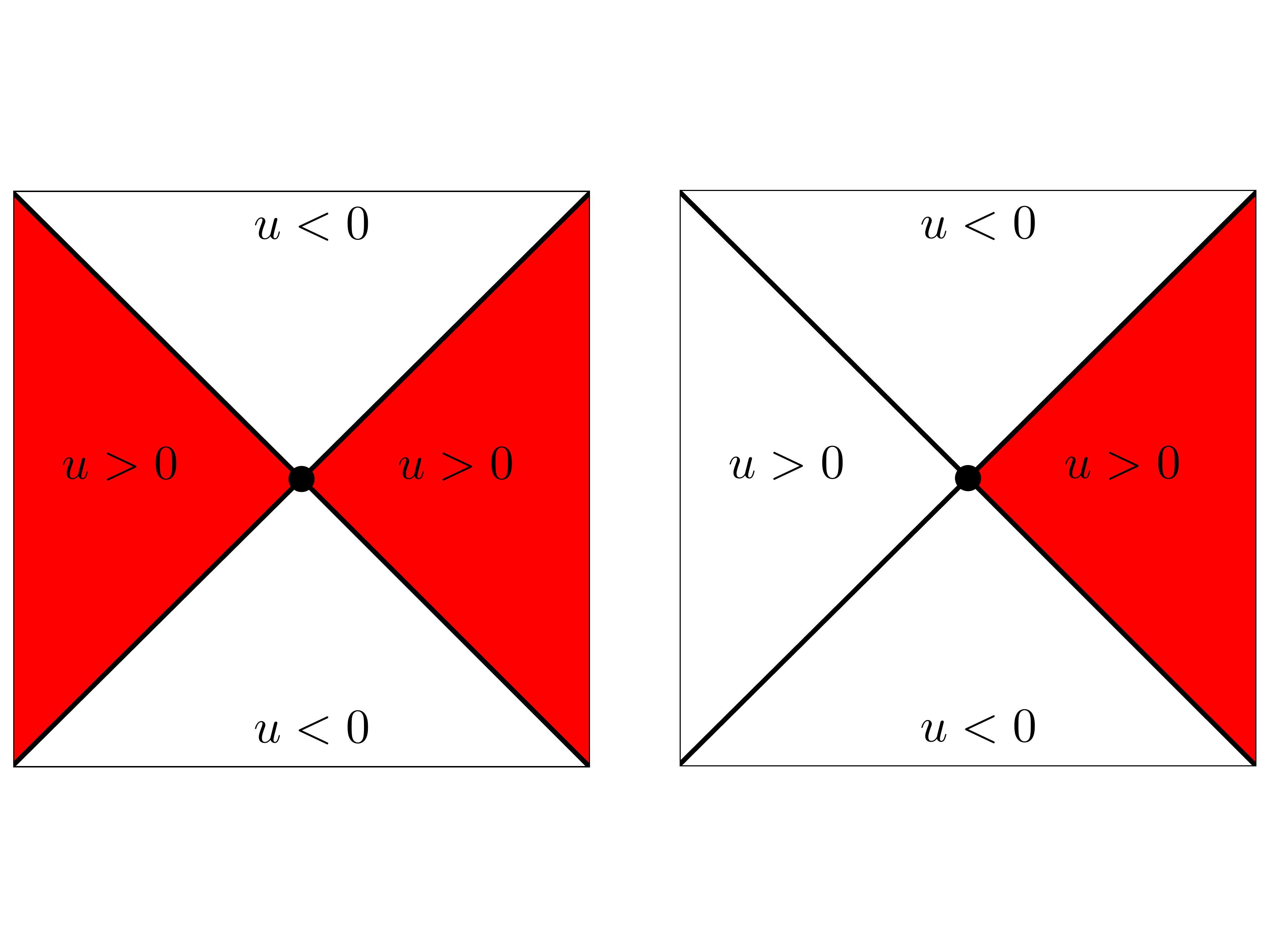}
\end{center}
\caption{Wedges (with the set $\Omega$ painted red)}\label{figi2}
\end{figure}

3. There is $\lambda\in (0,1)$ such that 
$u(x^0+tx)/t^2\to {1\over 2}\left(\lambda x_1^2 + (1-\lambda)x_2^2\right)$ 
as $t\to 0$. The complement of $\{ u>0\}$ and that of $\Omega$ is the single point
$x^0$.
\end{thma}
\section{Proof of the Main Result}
Let $u$ be a solution of (\ref{solution}).
By $L^p$- and $C^\alpha$-estimates
$u\in W^{2,p}_{\rm loc}(\R^3)\cap C^{1,\alpha}_{\rm loc}(\R^3)$ for each $p\in (1,+\infty)$ and $\alpha\in (0,1)$.
Differentiating $u$ we obtain
\begin{equation}\label{deriv}
\Delta \partial_k u = -\Phi'(u) \chi_\Omega \partial_k u\in L^\infty_{\rm loc}(\R^3)
\end{equation}
in the sense of distributions.
Consequently $u\in W^{3,p}_{\rm loc}(\R^3)\cap C^{2,\alpha}_{\rm loc}(\R^3)$ for each $p\in (1,+\infty)$ and $\alpha\in (0,1)$.
From (\ref{solution}) we infer now that the Hessian of $u$
satisfies
$$ |D^2 u| \ge c(n) >0 \textrm{ on } \{ u=0\}.$$
As $\{ u=0\} \cap \{ \nabla u\ne 0\}$ is by the implicit function theorem
locally a $C^{2,\alpha}$-surface ---the regularity of the surface can be improved
to real analyticity by the methods in \cite{CF}---, we will focus on the singular set
$S=\{ u=0\}\cap \{ \nabla u= 0\}$.
From now on we will assume that $u$ is axisymmetric, that
is $u=u(r,x_3)$
and confine ourselves thus to a two-dimensional analysis.

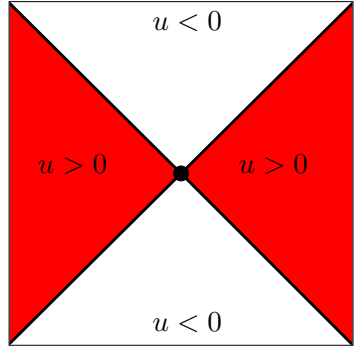
\begin{figure}
\begin{center}
\begin{picture}(0,0)%
\includegraphics{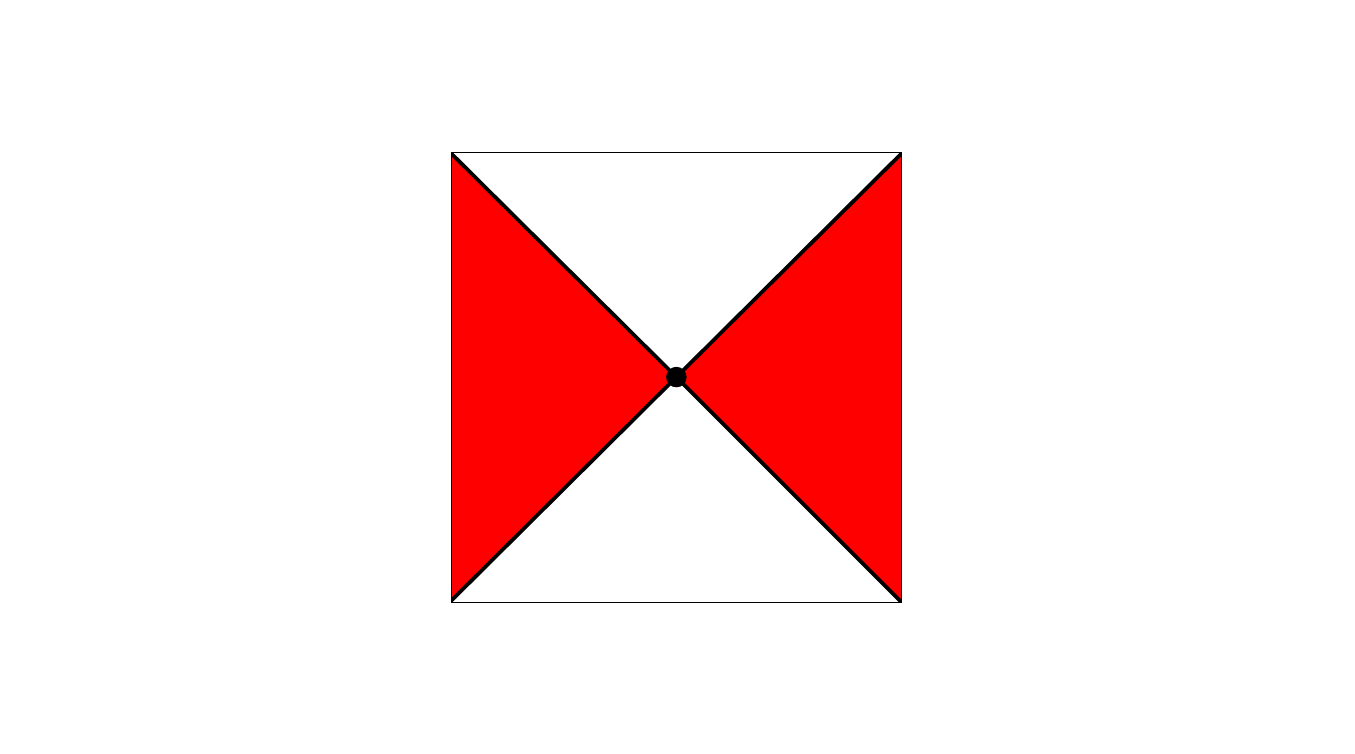}%
\end{picture}%
\setlength{\unitlength}{2368sp}%
\begingroup\makeatletter\ifx\SetFigFont\undefined%
\gdef\SetFigFont#1#2#3#4#5{%
  \reset@font\fontsize{#1}{#2pt}%
  \fontfamily{#3}\fontseries{#4}\fontshape{#5}%
  \selectfont}%
\fi\endgroup%
\begin{picture}(10824,6024)(-2411,-5173)
\put(2701,-3811){\makebox(0,0)[lb]{\smash{{\SetFigFont{11}{13.2}{\rmdefault}{\mddefault}{\updefault}{\color[rgb]{0,0,0}$u<0$}%
}}}}
\put(3601,-2161){\makebox(0,0)[lb]{\smash{{\SetFigFont{11}{13.2}{\rmdefault}{\mddefault}{\updefault}{\color[rgb]{0,0,0}$u>0$}%
}}}}
\put(1501,-2161){\makebox(0,0)[lb]{\smash{{\SetFigFont{11}{13.2}{\rmdefault}{\mddefault}{\updefault}{\color[rgb]{0,0,0}$u>0$}%
}}}}
\put(2701,-661){\makebox(0,0)[lb]{\smash{{\SetFigFont{11}{13.2}{\rmdefault}{\mddefault}{\updefault}{\color[rgb]{0,0,0}$u<0$}%
}}}}
\end{picture}%
\end{center}
\caption{Double Wedge (with the set $\Omega$ painted red)}\label{fig4}
\end{figure}

\begin{figure}
\begin{center}
\begin{picture}(0,0)%
\includegraphics{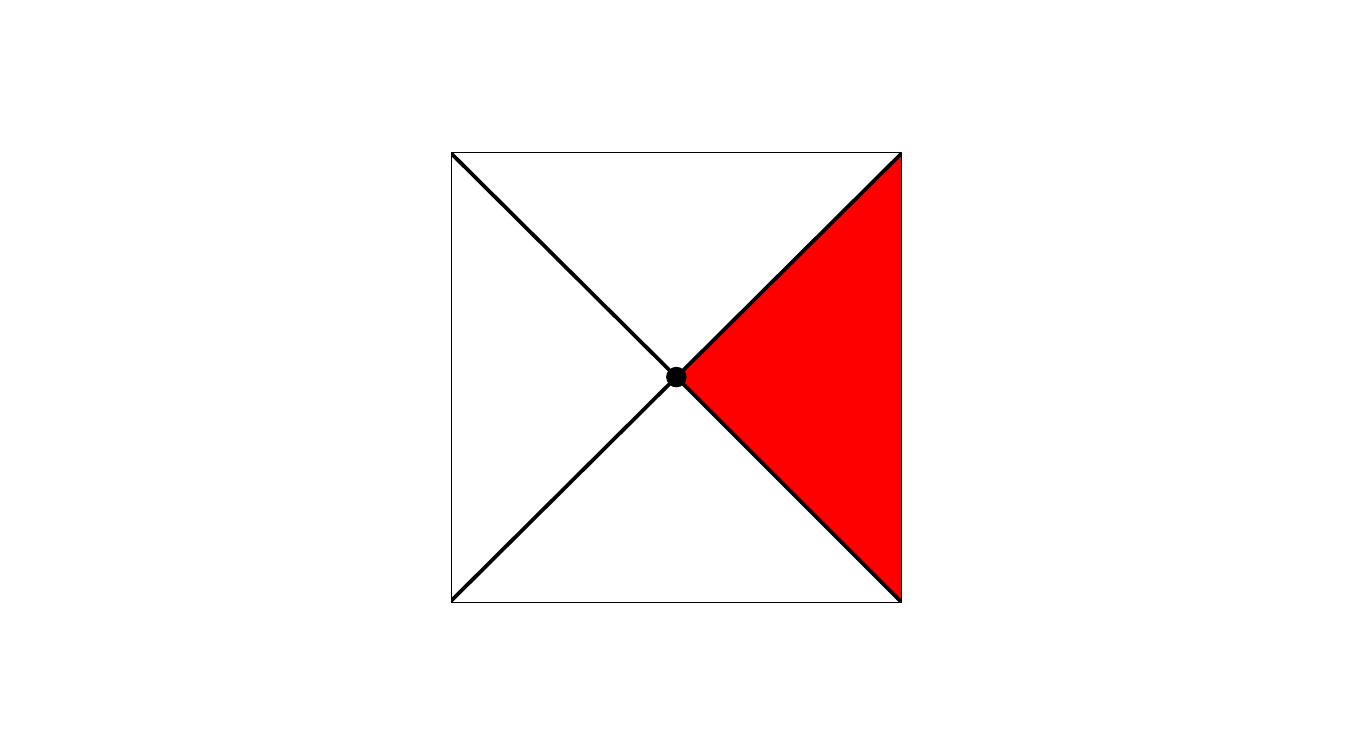}%
\end{picture}%
\setlength{\unitlength}{2368sp}%
\begingroup\makeatletter\ifx\SetFigFont\undefined%
\gdef\SetFigFont#1#2#3#4#5{%
  \reset@font\fontsize{#1}{#2pt}%
  \fontfamily{#3}\fontseries{#4}\fontshape{#5}%
  \selectfont}%
\fi\endgroup%
\begin{picture}(10824,6024)(-2411,-5173)
\put(2701,-3811){\makebox(0,0)[lb]{\smash{{\SetFigFont{11}{13.2}{\rmdefault}{\mddefault}{\updefault}{\color[rgb]{0,0,0}$u<0$}%
}}}}
\put(3601,-2161){\makebox(0,0)[lb]{\smash{{\SetFigFont{11}{13.2}{\rmdefault}{\mddefault}{\updefault}{\color[rgb]{0,0,0}$u>0$}%
}}}}
\put(1501,-2161){\makebox(0,0)[lb]{\smash{{\SetFigFont{11}{13.2}{\rmdefault}{\mddefault}{\updefault}{\color[rgb]{0,0,0}$u>0$}%
}}}}
\put(2701,-661){\makebox(0,0)[lb]{\smash{{\SetFigFont{11}{13.2}{\rmdefault}{\mddefault}{\updefault}{\color[rgb]{0,0,0}$u<0$}%
}}}}
\end{picture}%
\end{center}
\caption{Pseudo Wedge (with the set $\Omega$ painted red)}\label{fig5}
\end{figure}
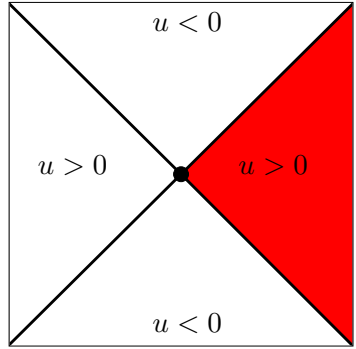

At each $x^0\in S$ we may rotate axes such that
$$D^2 u(x^0)=\begin{pmatrix}
\lambda&0\\
0&1-\lambda
\end{pmatrix}
\textrm{ for some }\lambda \in \R.$$
\noindent{\bf Case 1:}
If $0<\lambda<1$, then $\Omega^c$ consists in a sufficiently small ball $B_\delta(x^0)$
of only the point $x^0$ which is in this case a local minimum point of $u$.  
\\
{\bf Case 2:}
If $\lambda>1$ or $\lambda<0$,
$\{ u=0\}$ consists of two $C^1$-curves intersecting at a nonzero
angle at $x^0$ (cf. Figure \ref{fig4} and Figure \ref{fig5}): we may assume that $x^0=0$ and that $\lambda>1$.
As in this case for sufficiently small $\delta$,
$\partial_{x_1} u>0$ in $B_r\cap \{ x_1>\delta\}$ and
$\partial_{x_1} u<0$ in $B_r\cap\{x_1<-\delta\}$,
we may rescale and
obtain that $\{ u=0\}\setminus \{ 0\}$ consists of four $C^{1,\alpha}$-graphs.
The fact that $$u(tx)/t^2 \to {1\over 2}\left(\lambda x_1^2 + (1-\lambda) x_2^2\right)\textrm{ as } t \to 0$$
implies now
that the graphs have tangents as $x\to 0$ and that we may combine them to
two $C^1$-curves intersecting at a nonzero
angle at $0$.
\\
{\bf Case 3:}
If $\lambda=1$ or $\lambda=0$, then 
$\{ u=0\}$ consists either of two $C^1$-curves
ending in a cusp at $x^0$ (cf. Figure \ref{fig2}) or intersecting in a double cusp
at $x^0$ (cf. Figure \ref{fig3}):
we may assume $x^0=0$ and $\lambda=1$, implying that
$u(tx)/t^2 \to {1\over 2} x_1^2$ as $t\to 0$. From (\ref{deriv})
with $k=2$ we obtain that
$$ \Delta \partial_2 u = c(x) \partial_2 u$$
with H\"older continuous coefficients $c(x)$.
Applying \cite[Lemma 3.1]{fermi} repetitively 
for $\beta = 3/2, 7/2, 11/2, 15/2, \dots$
we infer
that either
$$\partial_2 u = p + \Gamma$$
where $p$ is a nontrivial harmonic polynomial
of degree $[\beta]+2$ 
with leading term of order $\ge 2$ (by the fact that
$u(tx)/t^2 \to {1\over 2} x_1^2$ as $t\to 0$)
and 
$$|\Gamma(x)|\le C_1 |x|^{\beta+2},$$
or $\partial_2 u$ vanishes of infinite order at $0$, that is
$$|\partial_2 u(x)| \le C_k |x|^k \textrm{ in }B_{r_k}(0)$$
for every $k\in \N$.
In the latter case we obtain by repetitive application of a 
well-known {\em strong unique continuation
property} (see \cite[Remark 6.7]{jk} for a very general result), that
$$\partial_2 u \equiv 0\textrm{ in each connected component of } \{ u\ne 0\} \textrm{ touching the
origin,}$$
implying by our information on the blow-up limit that
$$u \equiv {1\over 2}x_1^2 + f(x_1) \textrm{ in each connected component of } \{ u\ne 0\} \textrm{ touching the
origin,}$$
where $f(z)\to 0$ as $z\to 0$.
But this contradicts the constraint $\int_{\R^3}\rho=1$, thus proving that
infinite order vanishing is not possible.

Let us return to the former case
$$\partial_2 u = q + O(|x|^{k+{1\over 2}})$$
where $q$ is a nontrivial homogeneous harmonic polynomial
of degree $k\ge 2$.
It is important to note that ---by the fact that $q$ is harmonic--- if $q((0,1))=0$ then
$\partial_{x_1}q((0,1))\ne 0$. Similarly, if $q((0,-1))=0$ then
$\partial_{x_1}q((0,-1))\ne 0$.

It follows that
$$u={1\over 2}x_1^2+ \int_0^{x_2} q(x_1,s)\> ds  + O(|x|^{k+{3\over 2}}).$$
In this case
$\{ u<0\}$ is in a neighborhood of $0$ a one-sided or two-sided
cusp (depending on the signs of $q((0,1))$ and $q((0,-1))$); see Figure \ref{fig2}
and Figure \ref{fig3}.
In the special case that $q((0,1))=0$ or $q((0,-1))=0$,
we obtain $\partial_{x_1}q((0,1))\ne 0$ or $\partial_{x_1}q((0,-1))\ne 0$, respectively,
in which case we obtain a non-symmetric cusp (even asymptotically)
with respect to the $x_1$-axis
on the side of $(0,1)$ or $(0,-1)$, respectively.
\bibliographystyle{plain.bst}
\bibliography{stars.bib}

\begin{thebibliography}{10}

\bibitem{ACF}
Hans~Wilhelm Alt, Luis~A. Caffarelli, and Avner Friedman.
\newblock Variational problems with two phases and their free boundaries.
\newblock {\em Trans. Amer. Math. Soc.}, 282(2):431--461, 1984.

\bibitem{AW}
J.~Andersson and G.~S. Weiss.
\newblock Cross-shaped and degenerate singularities in an unstable elliptic
  free boundary problem.
\newblock {\em J. Differential Equations}, 228(2):633--640, 2006.

\bibitem{ASW2}
John Andersson, Henrik Shahgholian, and Georg~S. Weiss.
\newblock On the singularities of a free boundary through fourier expansion.
\newblock {\em Invent. Math.}

\bibitem{ASW}
John Andersson, Henrik Shahgholian, and Georg~S. Weiss.
\newblock Uniform regularity close to cross singularities in an unstable free
  boundary problem.
\newblock {\em Comm. Math. Phys.}, 296(1):251--270, 2010.

\bibitem{beals}
J.~F.~G. Auchmuty and Richard Beals.
\newblock Variational solutions of some nonlinear free boundary problems.
\newblock {\em Arch. Rational Mech. Anal.}, 43:255--271, 1971.

\bibitem{fermi}
Luis~A. Caffarelli and Avner Friedman.
\newblock The free boundary in the {T}homas-{F}ermi atomic model.
\newblock {\em J. Differential Equations}, 32(3):335--356, 1979.

\bibitem{CF}
Luis~A. Caffarelli and Avner Friedman.
\newblock The shape of axisymmetric rotating fluid.
\newblock {\em J. Funct. Anal.}, 35(1):109--142, 1980.

\bibitem{partial}
Luis~A. Caffarelli and Avner Friedman.
\newblock Partial regularity of the zero-set of solutions of linear and
  superlinear elliptic equations.
\newblock {\em J. Differential Equations}, 60(3):420--433, 1985.

\bibitem{C}
S.~Chandrasekhar.
\newblock {\em An introduction to the study of stellar structure}.
\newblock Dover Publications Inc., New York, N. Y., 1957.

\bibitem{C2}
S.~Chandrasekhar.
\newblock {\em Ellipsoidal figures of equilibrium}.
\newblock Dover Publications Inc., New York, N. Y., 1992.

\bibitem{CL}
Sagun Chanillo and Yan~Yan Li.
\newblock On diameters of uniformly rotating stars.
\newblock {\em Comm. Math. Phys.}, 166(2):417--430, 1994.

\bibitem{jk}
David Jerison and Carlos~E. Kenig.
\newblock Unique continuation and absence of positive eigenvalues for
  {S}chr\"odinger operators.
\newblock {\em Ann. of Math. (2)}, 121(3):463--494, 1985.
\newblock With an appendix by E. M. Stein.

\bibitem{libao}
Haigang Li and Jiguang Bao.
\newblock Existence of rotating stars with prescribed angular velocity law.
\newblock {\em Houston J. Math.}, 37(1):297--309, 2011.

\bibitem{yanyan}
Yan~Yan Li.
\newblock On uniformly rotating stars.
\newblock {\em Arch. Rational Mech. Anal.}, 115(4):367--393, 1991.

\bibitem{luo}
Tao Luo and Joel Smoller.
\newblock Rotating fluids with self-gravitation in bounded domains.
\newblock {\em Arch. Ration. Mech. Anal.}, 173(3):345--377, 2004.

\bibitem{MW}
R.~Monneau and G.~S. Weiss.
\newblock An unstable elliptic free boundary problem arising in solid
  combustion.
\newblock {\em Duke Math. J.}, 136(2):321--341, 2007.

\bibitem{tass}
J.L. Tassoul.
\newblock {\em Theory of rotating stars}.
\newblock Princeton Univ. Press, New Jersey, 1978.

\end{thebibliography}
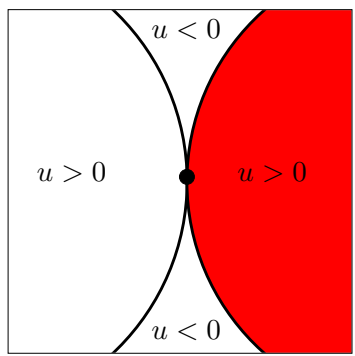
\begin{figure}
\begin{center}
\begin{picture}(0,0)%
\includegraphics{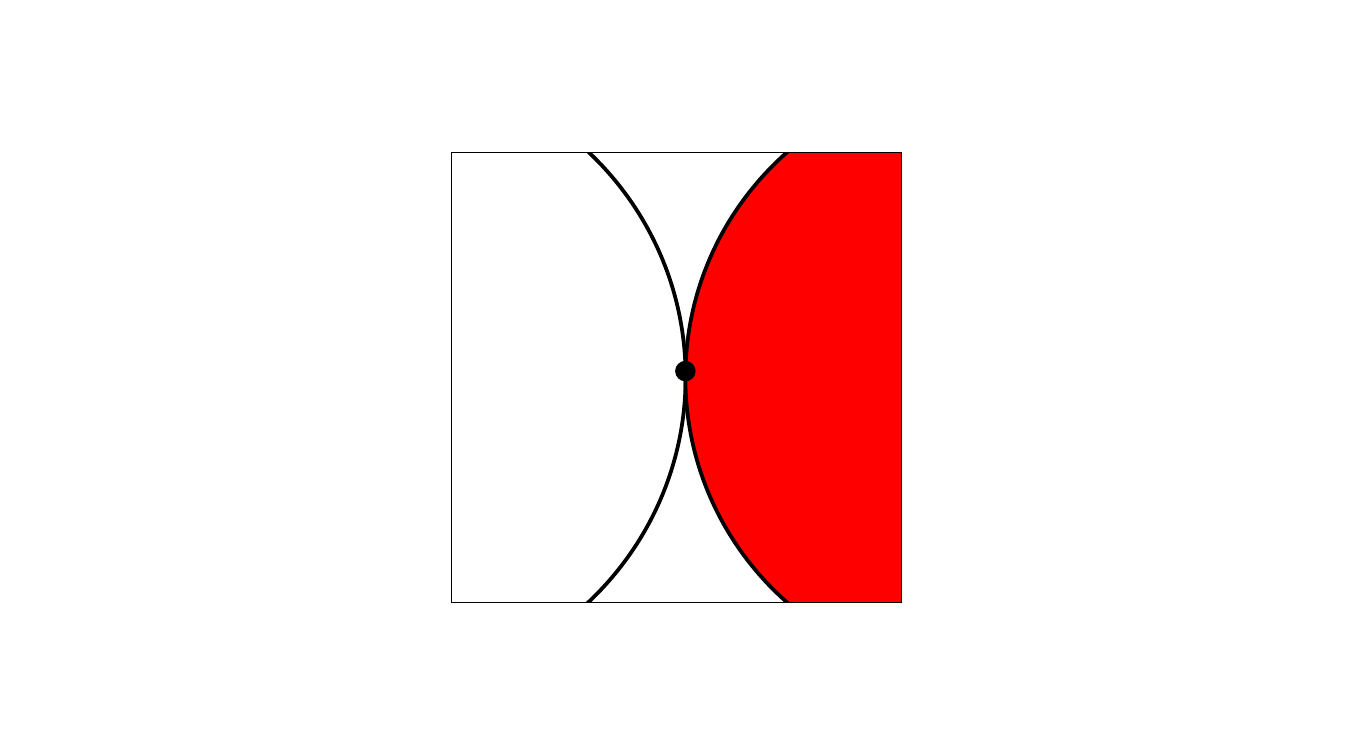}%
\end{picture}%
\setlength{\unitlength}{2368sp}%
\begingroup\makeatletter\ifx\SetFigFont\undefined%
\gdef\SetFigFont#1#2#3#4#5{%
  \reset@font\fontsize{#1}{#2pt}%
  \fontfamily{#3}\fontseries{#4}\fontshape{#5}%
  \selectfont}%
\fi\endgroup%
\begin{picture}(10824,6024)(-2411,-5173)
\put(2701,-3811){\makebox(0,0)[lb]{\smash{{\SetFigFont{11}{13.2}{\rmdefault}{\mddefault}{\updefault}{\color[rgb]{0,0,0}$u<0$}%
}}}}
\put(3601,-2161){\makebox(0,0)[lb]{\smash{{\SetFigFont{11}{13.2}{\rmdefault}{\mddefault}{\updefault}{\color[rgb]{0,0,0}$u>0$}%
}}}}
\put(1501,-2161){\makebox(0,0)[lb]{\smash{{\SetFigFont{11}{13.2}{\rmdefault}{\mddefault}{\updefault}{\color[rgb]{0,0,0}$u>0$}%
}}}}
\put(2701,-661){\makebox(0,0)[lb]{\smash{{\SetFigFont{11}{13.2}{\rmdefault}{\mddefault}{\updefault}{\color[rgb]{0,0,0}$u<0$}%
}}}}
\end{picture}%
\end{center}
\caption{Pseudo Cusp (with the set $\Omega$ painted red)}\label{fig1}
\end{figure}
\begin{figure}
\begin{center}
\begin{picture}(0,0)%
\includegraphics{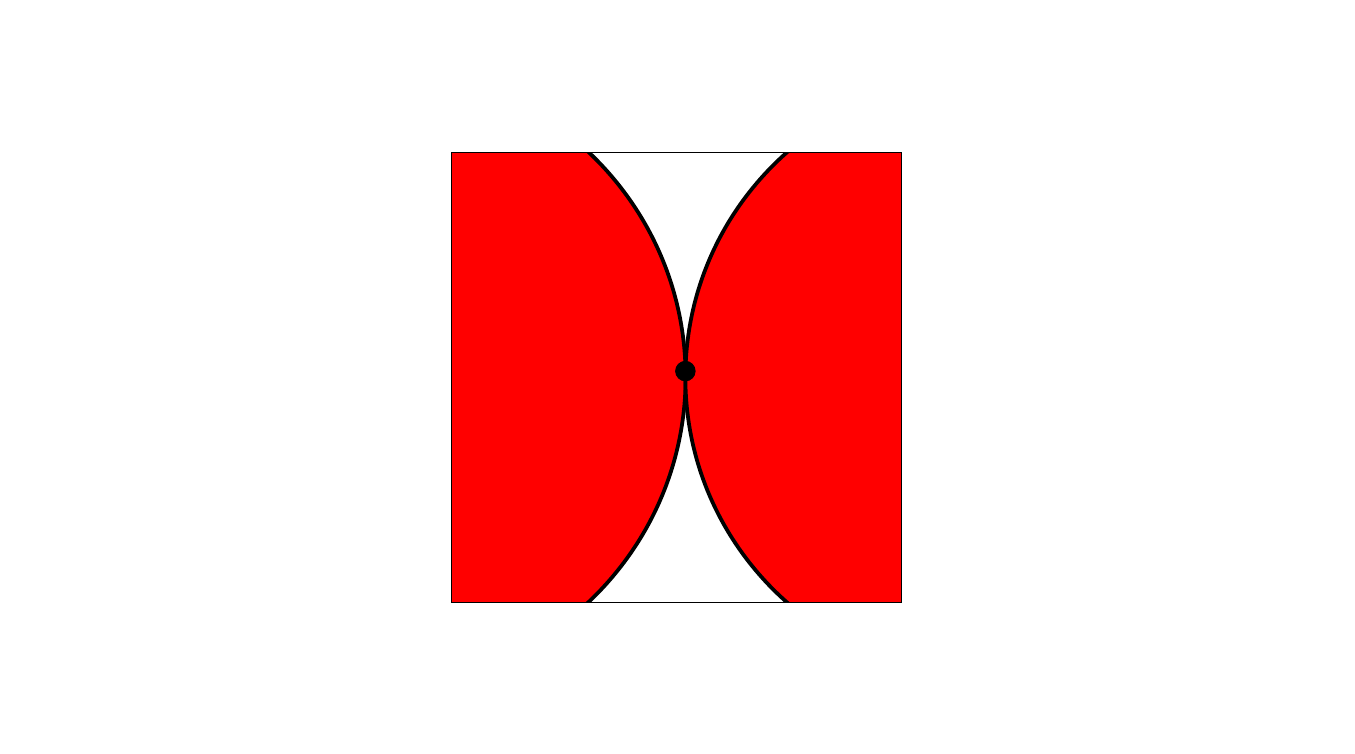}%
\end{picture}%
\setlength{\unitlength}{2368sp}%
\begingroup\makeatletter\ifx\SetFigFont\undefined%
\gdef\SetFigFont#1#2#3#4#5{%
  \reset@font\fontsize{#1}{#2pt}%
  \fontfamily{#3}\fontseries{#4}\fontshape{#5}%
  \selectfont}%
\fi\endgroup%
\begin{picture}(10824,6024)(-2411,-5173)
\put(2701,-3811){\makebox(0,0)[lb]{\smash{{\SetFigFont{11}{13.2}{\rmdefault}{\mddefault}{\updefault}{\color[rgb]{0,0,0}$u<0$}%
}}}}
\put(3601,-2161){\makebox(0,0)[lb]{\smash{{\SetFigFont{11}{13.2}{\rmdefault}{\mddefault}{\updefault}{\color[rgb]{0,0,0}$u>0$}%
}}}}
\put(1501,-2161){\makebox(0,0)[lb]{\smash{{\SetFigFont{11}{13.2}{\rmdefault}{\mddefault}{\updefault}{\color[rgb]{0,0,0}$u>0$}%
}}}}
\put(2701,-661){\makebox(0,0)[lb]{\smash{{\SetFigFont{11}{13.2}{\rmdefault}{\mddefault}{\updefault}{\color[rgb]{0,0,0}$u<0$}%
}}}}
\end{picture}%
\end{center}
\caption{Double Cusp (with the set $\Omega$ painted red)}\label{fig2}
\end{figure}
\begin{figure}
\begin{center}
\begin{picture}(0,0)%
\includegraphics{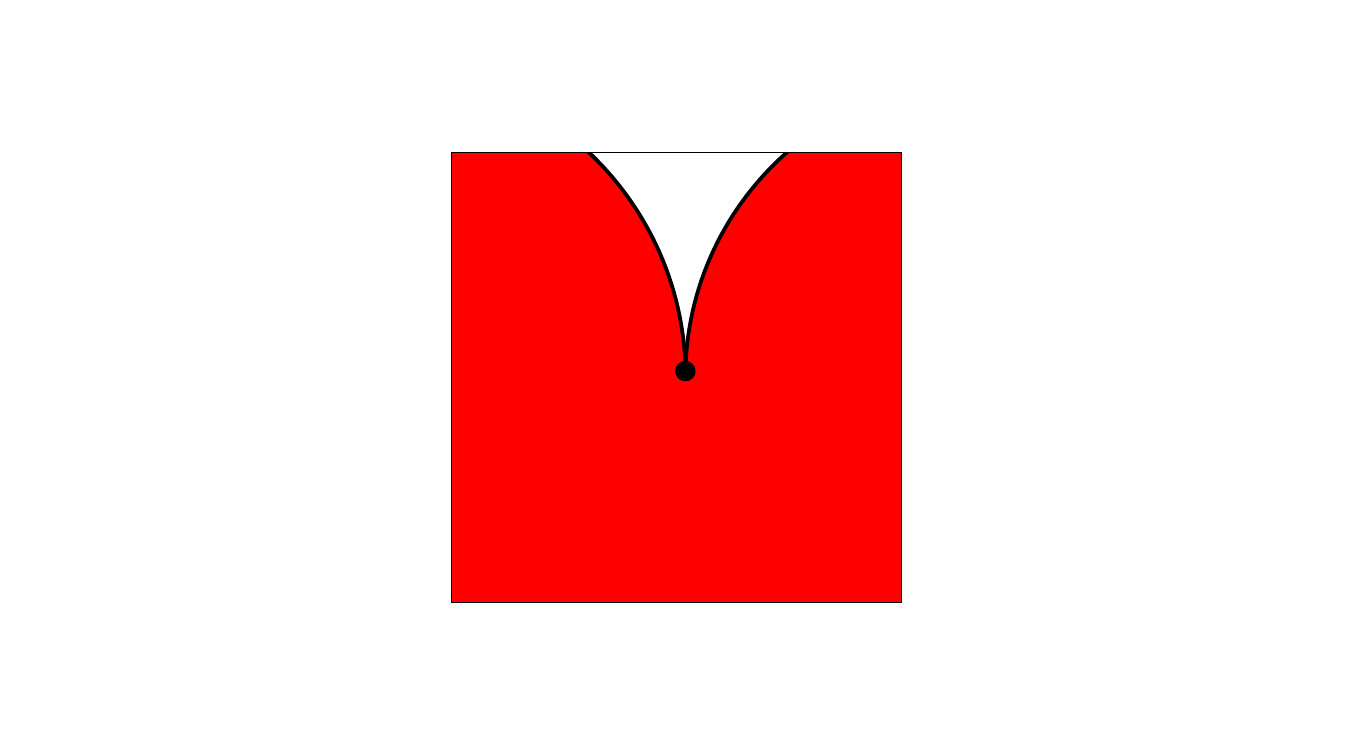}%
\end{picture}%
\setlength{\unitlength}{2368sp}%
\begingroup\makeatletter\ifx\SetFigFont\undefined%
\gdef\SetFigFont#1#2#3#4#5{%
  \reset@font\fontsize{#1}{#2pt}%
  \fontfamily{#3}\fontseries{#4}\fontshape{#5}%
  \selectfont}%
\fi\endgroup%
\begin{picture}(10824,6024)(-2411,-5173)
\put(2701,-3211){\makebox(0,0)[lb]{\smash{{\SetFigFont{11}{13.2}{\rmdefault}{\mddefault}{\updefault}{\color[rgb]{0,0,0}$u>0$}%
}}}}
\put(2701,-661){\makebox(0,0)[lb]{\smash{{\SetFigFont{11}{13.2}{\rmdefault}{\mddefault}{\updefault}{\color[rgb]{0,0,0}$u<0$}%
}}}}
\end{picture}%
\end{center}
\caption{Cusp (with the set $\Omega$ painted red)}\label{fig3}
\end{figure}
\end{document}